\documentclass[12pt]{article}
\usepackage{amssymb}
\usepackage{graphicx}
\usepackage{latexsym}
\def\part#1{\frac{\partial\phantom{q}}{\partial#1}}

\newenvironment{rmk}{\begin{trivlist}\item[]{\bf Remark:} }
{\end{trivlist}}

\newenvironment{ex}{\begin{trivlist}\item[]{\bf Example:} }
{\end{trivlist}}
\newenvironment{prf}{\begin{trivlist}\item[]{\bf Proof:} }
{\hfill $\Box$ \end{trivlist}}

\newtheorem{thm}{Theorem}
\newtheorem{definition}{Definition}
\newtheorem{prp}[thm]{Proposition}

\newcommand{\lie}[1]{\mathfrak{#1}}
\def\End{\mathop{\rm End}\nolimits}

\def\Hom{\mathop{\rm Hom}\nolimits}

\def\ker{\mathop{\rm ker}\nolimits}
\def\coker{\mathop{\rm coker}\nolimits}

\def\im{\mathop{\rm Im}\nolimits}
\def\deg{\mathop{\rm deg}\nolimits}
\def\rk{\mathop{\rm rk}\nolimits}
\def\tr{\mathop{\rm tr}\nolimits}

\def\ad{\mathop{\rm ad}\nolimits}

\newcommand{\R}{\mathbf{R}}
\newcommand{\C}{\mathbf{C}}

\newcommand{\PP}{{\rm P}}

\begin{document}
\title{Remarks on Nahm's equations}
\author{Nigel Hitchin\\[5pt]}

\maketitle

\centerline{\it{Dedicated to Simon Donaldson on the occasion of his 60th birthday}}
\vskip .5cm
\noindent {\bf Abstract:} Nahm's equations are viewed in a more general context where they appear as a vector field on a moduli space of ${\mathcal O}(2)$-twisted Higgs bundles on the projective line. Zeros of this vector field correspond to torsion-free sheaves on a singular spectral curve which we translate  in terms of a smooth curve in three-dimensional projective space. We also show how generalizations of Nahm's equations are required when the spectral curve is non-reduced and deduce the existence of non-classical conserved quantities in this situation. 

\section{Introduction}
Nahm's equations are the reduction of the self-dual Yang-Mills equations from four to one dimension and have played an important role in many parts of geometry and physics -- from the original study of magnetic monopoles to a vast collection of  constructions of hyperk\"ahler metrics. In fact, since Euler's equations for a spinning top form the simplest example one could say that they entered the scene centuries ago. They are equations for a triple of $n\times n$ matrix -valued functions of $t$:
$$\frac{dT_1}{dt}=[T_2,T_3],\qquad \frac{dT_2}{dt}=[T_3,T_1],\qquad  \frac{dT_3}{dt}=[T_1,T_2].$$
In \cite{NH1} they also appeared naturally in the context of generalized complex structures. The moduli space of  {\it generalized holomorphic bundles} on the projective line   has an  action of a one-parameter group of the fundamental  
B-field symmetry which is omnipresent in generalized geometry. In the simplest case this turns out to be equivalent to evolution via Nahm's equations. It was a throwaway remark in  \cite{NH1}, but here we consider it in more detail and in particular look for fixed points in the moduli space. 

The equations are integrable in the sense that they correspond (up to conjugation) to a linear flow on the Jacobian of an algebraic curve, the spectral curve. Put like that  it seems as if there are no fixed points, but what happens is that they occur 
for  singular or reducible spectral curves and therefore have a different flavour  from the more standard treatment of Nahm's equations. We show, using the twistor theory that lies behind the geometry of monopoles,  that a fixed point  corresponds to a curve together with a line bundle in projective 3-space. 

The situation where the spectral curve is non-reduced involves further features, both for Nahm's equations and the fixed points.  We consider the  case where it is a {\it ribbon} and in doing so introduce new conserved quantities for Nahm's equations beyond the coefficients of the equation of the spectral curve.

\section{Co-Higgs bundles and Nahm's equations}
\subsection{Generalized geometry}
One of the basic features of generalized geometry is the extra symmetry beyond diffeomorphisms provided by the action of a closed 2-form, the B-field.  The concept of a generalized complex structure (as in \cite{NH0},\cite{MG}) includes an ordinary complex structure and then closed 2-forms of type $(1,1)$ preserve the generalized complex structure and hence transform naturally associated geometrical objects. The notion of a generalized holomorphic bundle on a generalized complex manifold becomes quite concrete for an ordinary complex structure (\cite{MG},\cite{NH1}):
\begin{definition} Let $M$ be a complex manifold. A generalized holomorphic bundle is a   holomorphic vector bundle $V$ together with a holomorphic section $\phi$ of $\End V\otimes T$ such that $\phi\wedge \phi=0$ as a section of $\End V\otimes \Lambda^2T$.
\end{definition}

 Replacing $T$ by $T^*$ gives us Higgs bundles, so these are also called co-Higgs bundles.  We shall retain this terminology since ``generalized" will be used in a different way later on. We adopt the differential-geometric approach to holomorphic bundles by considering a fixed $C^{\infty}$ vector bundle $V$ and a holomorphic structure $A$ defined by an operator $\bar\partial_A:\Omega^0(M,V)\rightarrow \Omega^{01}(M,V)$ with $\bar\partial_A^2=0$. 

If $B$ is a closed $(1,1)$-form  then the interior product of the matrix-valued vector field  $\phi$ with $B$ gives $i_{\phi}B\in \Omega^{01}(M,\End V)$ and the B-field transform is the new holomorphic structure defined by the $\bar\partial$-operator 
$$\bar\partial_B=\bar\partial+i_{\phi}B$$ on the same $C^{\infty}$ bundle $V$. The three conditions $\bar\partial B=0,\bar\partial \phi=0,\phi\wedge\phi=0$ show that $\bar\partial_B^2=0$ which is the integrability condition for the holomorphic structure. The last two show that $\phi$, which is unchanged, is holomorphic with respect to this new structure.  As shown in \cite{NH1} if $B=\bar\partial\theta$ then the pairs $(V,\bar\partial,\phi)$ and $(V,\bar\partial_B,\phi)$ are holomorphically equivalent. 
\subsection{Nahm's equations}
We shall consider co-Higgs bundles in the one-dimensional case of $\PP^1$, studied in some detail in \cite{Rayan}. In this case there is only a one-dimensional choice of Dolbeault cohomology class in $H^1(\PP^1,K)$ for $B$. Choose a generator $[\omega]$. Generically, if $c_1(V)=0$ the bundle $V$ will be a trivial rank $n$ bundle and then we can write
$$\phi=(\phi_0+\phi_1 z+\phi_2 z^2)\frac{d}{dz}$$
where the $\phi_i$ are constant $n\times n$ matrices.  

\begin{thm} Let $(V,\psi)$ be a rank $n$ co-Higgs bundle over $\PP^1$ with $V$ holomorphically trivial and $B$ a $(1,1)$-form whose integral is non-zero. Then if $t$ lies in a neighbourhood of $0\in \C$ over which the holomorphic structure $\bar\partial_{tB}$ is trivial, there is a $t$-dependent choice of trivialization in which $\psi$  is represented by  $\phi(t)$ and the components of $\phi(t)$ satisfy the equations 
$$\frac{d\phi_0}{dt}=-\frac{1}{2}[\phi_1,\phi_0],\qquad \frac{d\phi_1}{dt}=[\phi_0,\phi_2],\qquad \frac{d\phi_2}{dt}=\frac{1}{2}[\phi_1,\phi_2]$$
and $\phi(0)=\psi$.
\end{thm}

\begin{rmk} Writing $\phi_0=-(T_1+iT_2), \phi_1=-2iT_3,\phi_2=-(T_1-iT_2)$ gives Nahm's equations
$$\frac{dT_1}{dt}=[T_2,T_3],\qquad \frac{dT_2}{dt}=[T_3,T_1],\qquad  \frac{dT_3}{dt}=[T_1,T_2].$$
\end{rmk}
\begin{prf} Triviality of the bundle $V$ means the existence of a gauge transformation $g(t)$ such that 
\begin{equation}
ti_{\psi}B=g^{-1}\bar\partial g \label{B}
\end{equation}
Any two are related by $\tilde g=hg$ where for each $t$, $h(t)$ is a constant matrix. The Higgs field $\psi$ then defines  a $t$-dependent $\phi(t)$ by 
\begin{equation}
\psi=g^{-1}\phi g \label{phi}
\end{equation}
Differentiating (\ref{B})  with respect to $t$ gives
$$i_{\psi}B=-g^{-1}\dot g g^{-1}\bar\partial g +g^{-1}\bar\partial \dot g$$
or, conjugating by $g$,
$$i_{\phi}B=-\dot g g^{-1}\bar\partial g g^{-1}+\bar\partial \dot g g^{-1}=\bar\partial(\dot g g^{-1}).$$
For $B$ take the standard volume form 
$$\omega=\frac{dz d\bar z}{(1+z\bar z)^2}$$
then
$$\bar\partial(\dot g g^{-1})=i_{\phi}B=\frac{1}{(1+z\bar z)^2}(\phi_0+\phi_1 z+\phi_2 z^2) d\bar z.$$
Integrating to give a regular integral gives 
$$\dot g g^{-1}=\frac{-1}{z(1+z\bar z)}(\phi_0+\phi_1 z+\phi_2 z^2)+\frac{\phi_0}{z}+c(t)$$
for a choice of constant matrix $c$. Take $c=\phi_1/2$ and then 
\begin{equation}
\dot g g^{-1}=\frac{-1}{z(1+z\bar z)}(\phi_0+\phi_1 z+\phi_2 z^2)+\frac{\phi_0}{z}+\frac{\phi_1}{2}
\label{gdot}
\end{equation}

Differentiating (\ref{phi}) with respect to $t$ gives 
$$0=-g^{-1}\dot g g^{-1}\phi+g^{-1}\dot\phi g +g^{-1} \phi \dot g$$ or
$$\dot\phi=[\dot g g^{-1},\phi]$$
and substituting from (\ref{gdot}) we obtain
$$\dot\phi=\left [\frac{\phi_0}{z}+\frac{\phi_1}{2},\phi\right].$$
Equating coefficents of $z$ gives the result.
\end{prf}

\begin{rmk} The choice of $c$ gives the symmetrical form of Nahm's equations arising from their origin where the $T_i$ lie in a compact Lie algebra, and $\PP^1$ is endowed with the real structure $z\mapsto -1/\bar z$. Taking $c=0$ instead gives the equations
$$\frac{d\phi_0}{dt}=[\phi_0,\phi_1],\qquad \frac{d\phi_1}{dt}=[\phi_0,\phi_2],\qquad \frac{d\phi_2}{dt}=0$$

\end{rmk}

If $V$ has degree $k$ where $0<k<n$ then the generic splitting type of a holomorphic structure is 
$V= {\mathcal O}^{k}(1)\oplus{\mathcal O}^{n-k}$ and the rank $k$  subbundle is uniquely determined. The structure group then reduces to a parabolic subgroup, the subgroup of $GL(n,\C)$ preserving a $k$-dimensional subspace. Now the Higgs field has the form 
$$\phi=\pmatrix{ A & B\cr
                           C & D}$$
where, in the affine coordinate $z$, $A,B,C,D$ are matrix-valued polynomials of degree $2,3,1,2$ respectively. We can then write 
$$\phi=(\phi_0+\phi_1 z+\phi_2 z^2+\phi_3z^3)\frac{d}{dz}$$
where $\phi_2$ lies in the parabolic subalgebra and $\phi_3$ in its nilradical. 

Applying the B-field action as above gives an integral 
$$\dot g g^{-1}=\frac{-1}{z(1+z\bar z)}(\phi_0+\phi_1 z+\phi_2 z^2+\phi_3 z^3)+\frac{\phi_0}{z}$$
which is regular at the origin. But 
$$\dot g g^{-1}=\pmatrix{ \alpha & \beta\cr
                           \gamma & \delta}$$
where $\alpha,\delta$ are functions, $\beta$ is a $C^{\infty}$ section of ${\mathcal O}(1)$ and $\gamma$ of ${\mathcal O}(-1)$. A term of the form $z^k/z(1+z\bar z)$ extends smoothly to a section of ${\mathcal O}(m)$ if $k\le m+2$, so given the degrees of $A,B,C,D$ this is well-defined on $\PP^1$. With the constant $c=0$ the equations are: 

\begin{equation}
\frac{d\phi_0}{dt}=[\phi_0,\phi_1],\qquad \frac{d\phi_1}{dt}=[\phi_0,\phi_2], \qquad \frac{d\phi_2}{dt}=[\phi_0,\phi_3], \qquad \frac{d\phi_3}{dt}=0.
\label{odd}
\end{equation}

\section{Moduli spaces and the Nahm flow}
\subsection{Moduli spaces}
Just as in the case of Higgs bundles, one can introduce the notion of stability into our situation and construct moduli spaces \cite{Rayan}. A co-Higgs bundle $(V,\phi)$ on $\PP^1$ is stable if for any $\phi$-invariant holomorphic subbundle $U\subset V$, $\deg U/\rk U< \deg V/\rk V$. In the case of equality the pair is semi-stable. Since $\phi$-invariance implies that $U$ is also preserved by  $\bar\partial_B=\bar\partial+i_{\phi}B$ stability is clearly invariant under B-field transforms. The space of S-equivalence classes of co-Higgs bundles (where S-equivalence means replacing the Harder-Narasimhan filtration of a semistable bundle by its graded version) is a well-defined non-compact algeraic variety and, as with vector bundles themselves, when the degree and rank are coprime it is smooth. Moreover, as with Higgs bundles, the coefficients $a_k$ of the characteristic polynomial 
$$
\det(x-\phi)=x^n+a_1x^{n-1}+\dots + a_n$$
define a proper map to a vector space $$W=H^0(\PP^1,{\mathcal O}(2))\oplus H^0(\PP^1,{\mathcal O}(4))\oplus \dots \oplus H^0(\PP^1,{\mathcal O}(2n)).$$
The B-field action therefore defines a canonical holomorphic vector field on this moduli space and we shall call this more general action from now on the {\it Nahm flow}. 

\begin{ex}
Take $V$ to be of rank $2$ and degree $(-1)$. Then in \cite{Rayan} it is shown that the moduli space of stable co-Higgs bundles with $\tr\phi=0$ is the universal elliptic curve
$${\mathcal S}=\{(z,w,c_0,c_1,\dots, c_4):w^2=c_0+c_1z+\dots +c_4z^4\}.$$
More invariantly, ${\mathcal S}\subset {\mathcal O}(2)\times H^0(\PP^1,{\mathcal O}(4))$ is the divisor of $w^2-\pi^*q(z)$ where $w$ is the tautological section of $\pi^*{\mathcal O}(2)$ on the total space of $\pi:{\mathcal O}(2)\rightarrow \PP^1$ and $q=c_0+c_1z+\dots +c_4z^4$ is a section of ${\mathcal O}(4)$. From $\cite{Rayan}$ stability implies that $V\cong {\mathcal O}\oplus {\mathcal O}(-1)$ and so, as above, the Higgs field is of the form 
$$\phi=\pmatrix{ a & b\cr
                           c & d}$$
                           where in particular $c\in H^0(\PP^1,{\mathcal O}(1))$. If $c=0$ then ${\mathcal O}\subset V$ is invariant which contradicts stability, so $c$ has a unique zero $z_0$. Since $a$ is a section of ${\mathcal O}(2)$ there is a map from the moduli space ${\mathcal M}$ to ${\mathcal S}$ by setting $w=a(z_0), c(z)=\det\phi$ and this is in fact an isomorphism. 
                           
Considering the Nahm flow, the last equation in  (\ref{odd}) gives $\phi_3=const.$ and this, as the coefficent of $z^3$, is strictly upper triangular          so we may take it to be 
$$\phi_3=            \pmatrix{ 0 & 1\cr
                           0 & 0}$$     
From the other equations we obtain, with 
$c(z)=c_0+zc_1, a(z)=a_0+a_1z+a_2z^2$, 
$$\dot c_0=2(c_1a_0-a_1c_0), \qquad \dot c_1=-2a_2c_0$$ and since $z_0=-c_0/c_1$ this gives 
$$\dot z_0=-2(a_0+a_1z_0+a_2z_0^2)=-2a(z_0).$$
Thus, at the points where $(z,c_0,\dots, c_4)$ are local coordinates on ${\mathcal S}$ the vector field is 
$$w\frac{\partial}{\partial z}.$$
The parameter $z$ fails to be part of a coordinate system if $w=0$ in which case  $w$ is a coordinate and since $w^2=q(z)$ the vector field has the local form 
\begin{equation}
\frac{q'(z)}{2}\frac{\partial}{\partial w}.
\label{localform}
\end{equation}

\end{ex} 
\subsection{Fixed points of the Nahm flow}
 In the example above a zero of the vector field occurs where $w=0$ and   from (\ref{localform}) we then have $q'(z)=0$ which is when the elliptic curve $w^2=q(z)$  is singular. Note that it also vanishes if  $q\equiv 0$: this is where the Higgs field is nilpotent. 

To see this in more generality, we note that in the original generalized geometry formulation, we have the pair $(\bar\partial _{tB}, \psi)$ where the holomorphic structure is varying and $\psi$ is fixed and so clearly $\det(w-\psi)$ is constant, so that the vector field is always tangent to the fibres of the map ${\mathcal M}\rightarrow W$. This means that the curve in ${\mathcal O}(2)$ defined by the equation $\det(w-\phi)=0$, the spectral curve $S$, is fixed along the flow. In particular,   the coefficients of the characteristic polynomial are constants of integration of Nahm's equations.

A  naive treatment of the integrability of Nahm's equations as in \cite{NH2} assumes that the spectral curve $S\subset {\mathcal O}(2)$ is smooth. In this case the co-Higgs bundle $(V,\phi)$ is obtained from a line bundle $L$ on $S$ as the direct image $V=\pi_*L, \phi=\pi_*w$, where $w$ is again the tautological section of $\pi^*{\mathcal O}(2)$, and then   $L$ is the cokernel of $\phi-w: \pi^*V(-2)\rightarrow \pi^*V$. If $V$ is of rank $n$ then the genus of $S$ is $g=(n-1)^2$, its canonical bundle $K_S\cong \pi^*{\mathcal O}(2n-4)$. By Grothendieck-Riemann-Roch if $L$ has degree $d$ then $\deg V=d+n-n^2$, so the original Nahm equations require $d=n^2-n$. The bundle $V$ is then trivial if and only if $V(-1)$ has no sections which is when $L(-1)$ of degree $g-1$ does not lie on the theta-divisor of $S$.

The Nahm flow then consists of tensoring $L$ by the one-parameter group of line bundles $U_t=\exp (t w [\omega])\in H^1({\mathcal O}(2),{\mathcal O}^*)$ restricted to $S$. Here $\omega$ is the standard $(1,1)$-form used in Theorem 1 and $[\omega]\in H^1(\PP^1,{\mathcal O}(-2))$ its cohomology class. The product with the tautological section $w$ of ${\mathcal O}(2)$ on its total space gives $w [\omega]\in H^1({\mathcal O}(2),{\mathcal O})$. Then $L\mapsto LU_t$ is a one-parameter group of translations in the Jacobian of $S$. Moreover, as in \cite{NH3}, the class $w[\omega]$ is always non-zero if  $\rk V>1$ and hence the flow has no fixed points. 

However, even in the original appearance of Nahm's equations for this author \cite{NH3}, singular and  reducible spectral curves  are  allowed, for example in the construction of axi-symmetric monopoles. Subsequent treatments of similar moduli spaces \cite{BNR}, \cite{CS},\cite{Sch}, identify the fibre as a compactified Jacobian parametrizing stable (in an appropriate sense) rank one torsion-free sheaves on the spectral curve. There is a large literature on compactified Jacobians but if we assume that the curve is reduced, then following \cite{Alex}, a torsion-free sheaf is given by the direct image of a line bundle on some partial normalization $S'$ of $S$. The generalized Jacobian $H^1(S',{\mathcal O}^*)$ of a singular curve is still a group so a fixed point of the Nahm flow must be represented by the direct image of a line bundle on a normalization $f:S'\rightarrow S$ for which the class $f^*w[\omega]=0\in H^1(S',{\mathcal O})$.

\begin{ex} In the example above the singular elliptic curves $w^2-q(z)=0$ are normalized by $\PP^1$ and $H^1(\PP^1,{\mathcal O})=0$ so any degree zero line bundle is trivial.
\end{ex}

Determining all such partial normalizations is seemingly a difficult task, but there is a more geometrical approach which we adopt now, and takes us back to the twistor theory of $\R^4$ and $\R^3$. 
\subsection{Twistor spaces and liftings}
Penrose's twistor theory encodes the Euclidean geometry of $\R^4$ in the holomorphic geometry of the complex 3-manifold ${\mathcal O}(1)\oplus {\mathcal O}(1)\rightarrow \PP^1$. The points of $\R^4$ correspond to holomorphic sections which are real with respect to an antiholomorphic involution with no fixed points. We are not concerned with reality here however. 

Any orientation-preserving Euclidean motion of $\R^4$ induces a holomorphic action on the twistor space, and in particular the one-parameter group   of translations $(x_0,x_1,x_2,x_3)\mapsto (x_0+t,x_1,x_2,x_3)$. The twistor space  is the complement of a line in $\PP^3$: in homogeneous coordinates $(z_0,z_1,z_2,z_3)$ we remove the line $z_0=z_1=0$ and then $[z_0,z_1]\in \PP^1$ defines the projection. The free holomorphic action is then
\begin{equation}
(z_0,z_1,z_2,z_3)\mapsto (z_0,z_1,z_2+t z_0,z_3-tz_1)
\label{action}
\end{equation}
and the invariant section $w=z_1z_2+z_3z_0$ of ${\mathcal O}(2)$ identifies the quotient by the action with the total space of ${\mathcal O}(2)$. As a principal $\C$-bundle over ${\mathcal O}(2)$ it defines a class  $\alpha\in H^1({\mathcal O}(2),{\mathcal O})$.

The quotient of $\R^4$ by the translation is $\R^3$ and each section of  ${\mathcal O}(1)\oplus {\mathcal O}(1)\rightarrow \PP^1$ projects to a section  $w=a_0z_0^2+a_1z_0z_1+a_2z_1^2=a(z_0,z_1)$ of ${\mathcal O}(2)\rightarrow \PP^1$. The three-dimensional space of such real sections is the twistor interpretation of the Euclidean geometry of $\R^3$ as in \cite{NH4}. A fixed section $w=a(z_0,z_1)$ of ${\mathcal O}(2)\rightarrow \PP^1$ has a one-parameter family of inverse images in  ${\mathcal O}(1)\oplus {\mathcal O}(1)$ and these sweep out a surface $$z_1z_2+z_3z_0-a(z_0,z_1)=0.$$
Adding in the line $z_0=z_1=0$ gives a smooth projective quadric in $\PP^3$ and the inverse images form one of the two families of lines. 

\begin{rmk} According to \cite{NH4}, a holomorphic vector bundle on ${\mathcal O}(2)$ trivial on each real section corresponds to a solution to the  Bogomolny equations $F_A=\ast \nabla \phi$ on $\R^3$. The class  $\exp \alpha\in H^1({\mathcal O}(2),{\mathcal O}^*)$ gives $A=0,\phi= 1$.
\end{rmk}

\begin{prp}The class $\alpha$ is (up to a multiple)  the same as the class $w[\omega]$ which gives the Nahm flow. 
\end{prp}
\begin{prf} Note that over the open set $U_0$ where $z_0\ne 0$ we have a section of the action defined by $(z_0,z_1, 0, z_3)$ and similarly over $U_1$ where $z_1\ne 0$ we have $(z_0,z_1, z_2, 0)$. Thus  a \v{C}ech cocycle in $H^1({\mathcal O}(2),{\mathcal O})$ defining it is provided by the value of $t$ on $U_0\cap U_1$ which relates these two sections. This is $t=w/z_0z_1$. 

Now $z_0z_1$ is the section of ${\mathcal O}(2)$ on $\PP^1$ vanishing at $0$ and $\infty$ and  using the affine parameter $z=z_1/z_0$, and identifying ${\mathcal O}(2)$ with the tangent bundle this is the vector field $zd/dz$. But its inverse, the form $dz/z$, is a  cocycle on $U_0\cap U_1$ which is a generator of $H^1(\PP^1,K)$, so the class $w[\omega]$ is represented by $w/z_0z_1$. 
\end{prf}
It follows that  if $C$ is a partial normalization of $S$ on which the pull-back of the class $w[\omega]$ is zero, then a choice of trivialization lifts it to a map into the principal $\C$-bundle over $S$. This is then a curve in $\PP^3$
which misses the line $z_0=z_1=0$. Conversely any such curve  projects to a curve $S$ in ${\mathcal O}(2)$ and  
points in $C$ which lie in the same orbit of the $\C$-action map to  singular points of $S$. This way $C$ is a partial normalization of $S$  and by construction  
 the class $w[\omega]$ is  trivial  on $C$. The direct image of any line bundle on $C$ is a torsion-free sheaf on $S$, and  taking the direct image on $\PP^1$ we have a  rank $n$ co-Higgs bundle where $n=\deg C$, whose equivalence class in the moduli space is fixed by the Nahm flow.  
\subsection{Commuting pairs}
In the generic case where the bundle $V$ on $\PP^1$ is trivial, the Nahm flow yields Nahm's equations
$$\frac{dT_1}{dt}=[T_2,T_3],\qquad \frac{dT_2}{dt}=[T_3,T_1],\qquad  \frac{dT_3}{dt}=[T_1,T_2].$$
and a zero of the induced vector field in the moduli space consists of matrices $(T_1,T_2,T_3)$ where a fourth matrix $T_0$ satisfies 
$$[T_0,T_1]=[T_2,T_3],\qquad [T_0,T_2]=[T_3,T_1],\qquad [T_0,T_3]=[T_1,T_2].$$
\begin{rmk} If  $T_0, T_1, T_2, T_3$  lie in the Lie algebra of a compact Lie group $G$ with a bi-invariant metric  then these equations are equivalent to the vanishing of the hyperk\"ahler moment map $\mu:\lie{g}\otimes {\mathbf H} \rightarrow \lie{g}\otimes \R^3$ for the adjoint action of $G$ on the flat hyperk\"ahler manifold $\lie{g}\otimes {\mathbf H}$. However we are dealing here with the complex case -- there are no non-trivial solutions for a compact group. We can see this by interpreting the equations as giving a translation-invariant solution to the self-dual Yang-Mills equations on $\R^4$, or equivalently a translation-invariant solution to the Bogomolny equations $\ast d_A\phi=F_A$ on $\R^3$. Quotienting by a lattice in $\R^3$ we have a solution on the 3-torus, but the Bianchi identity gives $0=d_AF_A=d_A\ast d_A\phi$. Integrating $(d_A\ast d_A\phi,\phi)$ and using Stokes' theorem we get $d_A\phi=F_A=0$.
\end{rmk}
To link this up with the above spectral curve approach we collect $T_1,T_2,T_3$, as in Theorem 1, into a co-Higgs field $\phi=\phi_0+\phi_1z+\phi_2 z^2$ and obtain (with $\psi=T_0$) 
$$\left[-\psi+\frac{\phi_0}{z}+\frac{\phi_1}{2},\phi\right]=0.$$
Hence the  term $\phi_-=-z\psi+{\phi_0}+z{\phi_1}/{2}$, which is  linear in $z$, defines a matrix  with entries in $H^0(\PP^1,{\mathcal O}(1))$ which commutes with $\phi$. Consider also
$$\phi_+=\frac{\phi_1}{2}+\phi_2 z+\psi$$
which is a similar section. This also  commutes with $\phi$ since $\phi_-+z\phi_+=\phi$.

Thus $(\phi_+,\phi_-)$ defines a matrix-valued section $\varphi$ of ${\mathcal O}(1)\oplus {\mathcal O}(1)$, and since $[\phi_+,\phi_-]=0$  we have $\varphi\wedge\varphi=0$ -- rather like a higher-dimensional Higgs field. Following this approach (due to Simpson \cite{CS} and in this context as in \cite{NH1}) it defines a  sheaf with compact support  on the total space of ${\mathcal O}(1)\oplus {\mathcal O}(1)$, or $\PP^3\backslash \PP^1$. If we denote by $x,y$ the tautological sections of ${\mathcal O}(1)$ on the two factors then $x$ acts by $\phi_+$, $y$ by $\phi_-$ and 
the  sheaf  is supported  on the variety defined by 
$$\det (u(x-\phi_+(z))+v(y-\phi_-(z)))=0$$
for all $u,v$. Roughly speaking it is the common cokernel of the family of commuting matrices $u(x-\phi_+(z))+v(y-\phi_-(z))$ and  is a rank $1$ sheaf supported  on the curve $C\subset \PP^3\backslash \PP^1$ above. Moreover, since $\phi_-+z\phi_+=\phi$ we have, putting $u=z,v=1$
$$\det(xz+y-\phi)=0$$
which with $w=xz+y$ is the equation of the spectral curve $S$ of $\phi$. This provides the projection to $S\subset {\mathcal O}(2)$. 
\subsection{Rank 2}\label{rk2}
 Consider the basic example where $\phi$ takes values  in $\lie{sl}(2,\C)$.  We take the equations for a fixed-point of the Nahm flow  in the form
$$[\psi,\phi_0]=\frac{1}{2}[\phi_0,\phi_1],\qquad [\psi,\phi_1]=[\phi_0,\phi_2],\qquad [\psi,\phi_2]=\frac{1}{2}[\phi_1,\phi_2].$$
Since $\phi$ becomes nilpotent at some point, without loss of generality we can take 
$$\phi_0=            \pmatrix{ 0 & 1\cr
                           0 & 0}.$$
                           The first equation gives $\psi+\phi_1/2=a\phi_0$. Substituting in the second we get $\phi_2-a\phi_1=b\phi_0$ and in the third $(a^2+b)[\phi_0,\phi_1]=0$ so either $[\phi_0,\phi_1]=0$ or $a^2+b=0$. In the first case, $\phi_1$ and $\phi_2$ are multiples of $\phi_0$ which means $\phi$ is nilpotent which we consider later.
                           
                            So with $a^2+b=0$ we have 
                           $\psi=a\phi_0-\phi_1/2$ and $\phi_2=a\phi_1-a^2\phi_0.$
                           This means 
                           $$\phi_-=(1-az)\phi_0+z\phi_1,\qquad \phi_+=a(1-az)\phi_0+az\phi_1$$
       so $\phi_+=a\phi_-$.            Moreover      $\tr \phi_-^2=    2z(1-az)\tr \phi_0\phi_1+z^2\tr\phi_1^2$.
       
       So the curve $C\subset \PP^3\backslash \PP^1$ has the equation in affine coordinates 
       $$y=ax,\qquad x^2=    2z(1-az)\tr \phi_0\phi_1+z^2\tr\phi_1^2.$$
       If $\tr \phi_0\phi_1\ne 0$, this is a nonsingular conic in the plane $y=ax$.
       
       Suppose $(x,y)$ and $(x+t,y-zt)$ lie on $C$. Then since $y=ax$, $z=-a$ and 
       $$x^2=-2a(1+a^2)\tr \phi_0\phi_1+a^2\tr\phi_1^2$$
 so in general there are two such points and the image $S$ has a double point where $w=xz+y=0$. If the right hand side is zero, then the vector field is tangential to $C$ and the image has a cusp. If $\tr \phi_0\phi_1= 0$ the curve $C$ is a pair of lines meeting in one point $x=y=z=0$. The image is a pair of sections of ${\mathcal O}(2)$ meeting at $(w,z)=(0,0)$ and $(w,z)=(0, -a)$. In all cases these are  partial normalizations with $H^1(C,{\mathcal O})=0$.

Now turn to the other zero of the vector field in the example above: where $\phi$ is nilpotent and the characteristic polynomial is $w^2$. The spectral curve in this case is the zero section of ${\mathcal O}(2)$ with multiplicity $2$: its first order neighbourhood. Let $X$ be the curve $w^2=0$ and $\PP^1=X_{red}$ the reduced curve, then there is an exact sequence of sheaves 
$$0 \rightarrow {\mathcal O}(-2)\rightarrow {\mathcal O}_X\rightarrow {\mathcal O}\rightarrow 0$$
and $H^1(X, {\mathcal O}_X)\cong H^1(\PP^1,{\mathcal O}(-2)) \cong \C$. In fact our class $w[\omega]$ is a generator.
\begin{rmk} In  the twistor theory  of monopoles  the Higgs field $\phi$ for a solution of the Bogomolny equations is precisely the obstruction to extending the trivialization of the corresponding holomorphic vector bundle on ${\mathcal O}(2)$ to the first order neighbourhood. In our case $\phi=1$ which is everywhere non-vanishing and  hence is a non-zero element of $H^1(X, {\mathcal O}_X)$. 
\end{rmk} 
We learn nothing more about the co-Higgs bundle from the spectral curve, but there is  extra information in the Higgs field $\phi=a(z)\phi_0$. In a neighbourhood of a point where $z\ne 0$ the cokernel of $\phi$ defines an invertible sheaf on $X$, generated by the cokernel of the constant matrix $\phi_0$. This is no longer true where $a(z)$ vanishes and all we get is a rank one torsion-free sheaf. Although the curve is smooth, we are in a similar situation to the general case and we can define 
$$\phi_{+}=(a_1/2+za_2)\phi_0,\qquad \phi_-=(a_0+za_1/2)\phi_0$$
where $\phi_-+z\phi_+=\phi$. Then  $\phi_-,\phi_+$ map $S$  into a curve $C\subset {\mathcal O}(1)\oplus {\mathcal O}(1)$ with equation $x^2=0=y^2$. Moreover if $a_1^2-4a_0a_2\ne 0$, $\phi_+,\phi_-$ have no common zero and the cokernels define a line bundle on $C$ whose direct image on $S$ is the required torsion-free sheaf. As above, the class $w[\omega]$ is trivial on $C$ and so we have a fixed point of the Nahm flow. When $a(z)$ has a double zero we take the direct image of a torsion-free sheaf, which is still invariant under tensoring by the line bundes. 
\section{Ribbons}
\subsection{Ribbons and line bundles}
The previous example is part of a more general picture where the spectral curve is non-reduced.  We restrict attention to multiplicity $2$ and a smooth reduced curve: this is called a {\it ribbon} \cite{BE}.
\begin{definition} A ribbon $X$ on  $S$ is a curve such that $X_{red}\cong S$ and the ideal sheaf ${\mathcal I}$ of $S$ in $X$ is invertible and satisfies ${\mathcal I}^2=0$.
\end{definition}
This is an abstract ribbon. We are concerned with a curve $X$ defined  by $\det(w-\phi)=p(w,z)^2=0$ in ${\mathcal O}(2)$ so that ${\mathcal I}$ is the conormal bundle of $S$ defined by $p(w,z)=0$. 

Simpson's results on the moduli spaces of sheaves imply \cite{CK} that if  $\det(w-\phi)=0$ defines a ribbon $X$ in the surface ${\mathcal O}(2)$, then the co-Higgs bundle is defined by the direct image of one of two types of sheaves:
\begin{itemize}
\item
a rank $2$ vector bundle $E$ on the reduced curve $X_{red}=S$
\item
a {\it generalized line bundle} on $X$,  a torsion-free sheaf which is  free of rank one outside a divisor $D\subset S$.
\end{itemize}
The first case is  rank one since ${\mathcal O}_S(E)$ and ${\mathcal O}_X$ have the same dimension as ${\mathcal O}_S$-modules, or equivalently the rank term in the Hilbert polynomial is $1$. In the second case  it was shown in \cite{BE} that there is a canonical blow up   $f:X'\rightarrow X$ of $X$ at the points of $D$ giving a ribbon $X'$ and the   generalized line bundle is then $f_*L$ for a line bundle $L$ on $X'$.

\begin{rmk} The first case occurs naturally in the Higgs bundle description of the moduli space of representations of a surface group into certain real Lie groups associated to the quaternions \cite{HS}.
\end{rmk}

\begin{ex} In Section \ref{rk2} the cokernel of $a(z)\phi_0$ defines a generalized line bundle on the ribbon $w^2=0$: the first order neighbourhood of the zero section $S$. There we lifted the curve to a quadric surface in $\PP^3$ where, if the two zeros of $a(z)$ were distinct, we had a line bundle. The blow-up in this case is achieved in the ambient surface ${\mathcal O}(2)$, which compactifies to the Hirzebruch surface $\PP({\mathcal O}\oplus {\mathcal O}(2))$. The zero section $S\subset \PP({\mathcal O}\oplus {\mathcal O}(2))$ has self-intersection $2$ and blowing up the two zeros of $a(z)$ this becomes zero. But the two $\PP^1$ fibres now have self-intersection $-1$ and can thus be blown down giving the quadric surface $\PP^1\times \PP^1$.
\end{ex}

We see from this that the Nahm flow is obtained by either tensoring the rank $2$ bundle $E$ on $S$ by $U_t$ or the line bundle $L$ by $f^*U_t$. When $S$ is smooth and has genus $> 0$ there are clearly no fixed points, and the example in Section \ref{rk2} shows what happens in the case of genus $0$.

\subsection{Ribbons and conserved quantities}
We shall discuss here the implications for Nahm's equations themselves when the spectral curve is a ribbon. The simplest case is where $V$ is the direct image of a rank $2$ bundle $E$ on $S$. This is when the Higgs field $\phi$  has 2-dimensional eigenspaces, or equivalently $p(\phi)=0$. The Nahm flow is then described by $E\mapsto E\otimes U_t$ and so the projective bundle $\PP(E)$ on $S$ is an invariant of the flow -- a geometric conserved quantity. 

The case of a generalized line bundle occurs when the generic eigenspaces are one-dimensional. To see what this means for the Nahm flow,  we follow the approach of Lucas Branco, who considers in his forthcoming Oxford DPhil thesis the Higgs bundle case.  

In our language we suppose then that we have a co-Higgs bundle $(V,\phi)$ on $\PP^1$ where $\rk V=2m$, $\det (w-\phi)=p^2(w)$ and $p=0$ defines a smooth curve $S$.  The two cases  of rank one sheaves on a ribbon correspond to whether the generic minimal polynomial is $p$ (the first case) or $p^2$. In the latter case, the kernel of $p(\phi)\in H^0(\PP^1,\End V(2m))$ defines a $\phi$-invariant subbundle $W_1\subset V$. Since $S$ is irreducible, there are no further invariant subbundles and since the generic minimal polynomial is of degree $m$ we have  $\rk W_1=m$. Thus  $V$ is an extension of co-Higgs bundles
$$0\rightarrow W_1\rightarrow V\rightarrow W_2\rightarrow 0$$
where $W_1,W_2$ have the same spectral curve $S$. With respect to a $C^{\infty}$-splitting we can therefore write
$$\bar\partial_V=\pmatrix{ \bar\partial_1  & \beta\cr
                           0 & \bar\partial_2}\qquad \phi=\pmatrix{ \varphi_1  & \psi\cr
                           0 & \varphi_2}$$
                           where $\bar\partial_V\phi=0$ implies 
                           \begin{equation}       
                          \bar\partial \psi+\varphi_{12}\beta=0.  \label{f12}
                      \end{equation}
                      In this equation for  $\beta\in \Omega^{01}(\PP^1, \Hom(W_2,W_1))$ we define $\varphi_{12}\beta=\beta\varphi_2-\varphi_1 \beta$.
                      
                           Since both $\bar\partial_V$ and $\phi$ preserve $W_1$ the   B-field action $\bar\partial_V\mapsto\bar\partial_V+i_{\phi}B$ also preserves the structure of an extension. 
                           
 \begin{rmk}         Since $W$ remains an extension, we have two Nahm flows corresponding to the bundles $W_1,W_2$. If $W$ is stable then $W_1$ has negative degree $-d$ and so $\deg W_2=d$. This means (unless $m$ divides $d$) that the Nahm flows on $W_1,W_2$ generically correspond to the equations (\ref{odd}) rather than the original Nahm equations, even if $W$ itself is holomorphically trivial.
  \end{rmk}
                          
                         Consider the holomorphic map $\varphi_{12}:\Hom(W_2,W_1)\rightarrow \Hom(W_2,W_1)(2)$. Its kernel and cokernel  are in fact holomorphic vector bundles. To see this note that $\varphi_1$ and $\varphi_2$ are Higgs fields with the same spectral curve $S$, so there are line bundles $L_1,L_2$ on $S$ whose direct images are $W_1,W_2$ and the Higgs fields are the direct images of $w:L_1\rightarrow L_1(2)$,  $w:L_2\rightarrow L_2(2)$. For a small open set $U\subset \PP^1$, $L_1$ and $L_2$ are isomorphic on $\pi^{-1}(U)$ and hence $W_1\cong W_2=W$ and  $\varphi_1=\varphi_2=\varphi$.  Trivializing ${\mathcal O}(2)$ over $U$,  $\ker \varphi_{12}$ can be identified with the sheaf of centralizers of a holomorphic matrix $\varphi$. Since $S$ is assumed smooth, $\varphi$ is regular and the space of centralizers is       spanned by $1,\varphi,  \dots,\varphi^{m-1}$. So globally $\ker \varphi_{12}$    is a rank $m$ holomorphic vector bundle and     the same holds for the cokernel.     
                          Equation (\ref{f12}) now says that the projection $\tilde\psi$ of $\psi$ to $\coker \varphi_{12}$ is holomorphic.
                          
                          Now the B-field action changes the holomorphic structure to $\bar\partial +i_{\phi}\omega$ and the induced operator on   $\psi\in \Omega^0(\PP^1,\Hom(W_2,W_1)(2))$ is 
                          $\bar\partial \psi +\varphi_{12}\omega \psi$. Since  $\varphi_{12}\omega\psi$ is trivial on the cokernel  the B-field action induces the same holomorphic structure on $\coker \varphi_{12}$.  In our formalism $\phi$ is unchanged and so the holomorphic section $\tilde\psi$ unchanged.   This  can therefore be considered as a conserved quantity under the Nahm flow.  (Strictly speaking $\psi$  is defined by the extension  rather than the bundle $V$ itself and so the invariant is the section up to a constant multiple). 
                          
                          \begin{rmk} A more sophisticated interpretation of the above is via the second spectral sequence of the hypercohomology for the complex of sheaves $\varphi_{12}:{\mathcal O}(\Hom(W_2,W_1))\rightarrow {\mathcal O}(\Hom(W_2,W_1)(2))$ \cite{LB}.
                          \end{rmk}
                          
      If  $\psi$ projects to zero in $\coker {\varphi}_{12}$ 
      then $\psi={\varphi}_{12}\theta$ for some $\theta\in \Omega^0(\PP^1,\Hom(W_2,W_1))$. But $\theta$ can  be used to change the $C^{\infty}$ splitting making $\psi=0$. In this case 
      $$\phi=\pmatrix{ \varphi_1  & 0\cr
                           0 & \varphi_2}$$

   and    $p(\phi)\equiv 0$ and we are back to the first case, so $\tilde\psi$, invariant by the flow, must be  part of the data of a generalized line bundle on $X$.  We shall see next what it is in a more concrete fashion next.

                          A local holomorphic section of $L_1^*L_2$ on $S$ defines a map from $L_1$ to $L_2$ commuting with the scalar multiplication by $w\in H^0(S,\pi^*({\mathcal O}(2)))$. The direct image therefore intertwines $\phi_1$ and $\phi_2$ and it follows  that $\ker \varphi_{12}\cong \pi_*(L_1^*L_2)$. Then $(\coker \varphi_{12})^*\cong \pi_*(L_1^*L_2)(-2)$. Relative duality gives
                          \begin{equation}
                          \coker \varphi_{12}\cong (\pi_*(L _1^*L_2))^*(2)\cong \pi_*(L_1L_2^*K_S)(4)\cong \pi_*(L_1L_2^*)(2m) \label{cok}
                          \end{equation}
    Thus the projection $\tilde\psi$  defines a non-zero holomorphic section $s$ of $L_1L_2^*(2m)$ on $S$, and this  vanishes on a divisor $D$, and since $\psi$ was really only defined up to a multiple it is the divisor which is the conserved quantity. Note that if $\ell_1,\ell_2$ are the degrees of the line bundles then $-d=\deg W_1=\ell_1+m-m^2, d=\deg W_2=\ell_2+m-m^2$ and so $\deg D=\ell_1-\ell_2+2m^2=-2d+2m^2$, so $0<d< m^2$.

     \begin{prp} Let $(w=\lambda, z=a)$ be a point of the divisor $D$ on the curve $S$. Then the $\lambda$-eigenspace of $\phi(a)$ has multiplicity $2$.
     \end{prp}
         \begin{prf} The application of relative duality in Equation \ref{cok} identifies $\coker \varphi_{12}$ with $\ker \varphi_{21}$ for a homomorphism $\varphi_{21}$ from $W_2$ to $W_1(2m)$. In fact, as we have seen, locally $\varphi_{12}$ can be considered as a holomorphic $m\times m$ matrix acting as $x\mapsto [\varphi,x]$. Using the invariant inner product $\tr(xy)$, $\ker \varphi$ is the orthogonal complement of the image of $\ad \varphi$. This maps to $\coker \varphi$ isomorphically unless the kernel of $(\ad \varphi)^2$ has dimension greater than $m$. But these points correspond to the discriminant locus of $p$, giving  the ramification points of $\pi:S\rightarrow \PP^1$, which provide the twist in  the relative duality formula.
         
       Pulling back $\pi^*(L_1L_2^*(2m))$ to $S$ there is the natural evaluation map 
        $$    \pi^*(L_1L_2^*(2m))_{(\lambda,a)}\rightarrow  L_1L_2^*(2m)_{(\lambda,a)}$$    
        and a  point $(\lambda, a)$ of $D$ is where the global section $s$ of $L_1L_2^*(2m)$ vanishes which means that at this point the direct image of $s$, $\tilde \psi$, maps the
        cokernel $L_2$ of $\phi_2-\lambda$ to zero in the cokernel $L_1$ of  $\varphi_1-\lambda$. Equivalently, $\im \psi\subseteq \im (\varphi_1-\lambda)$.
        
        Let $v_2$ be a $\lambda$-eigenvector of $\varphi_2$ at $z=a$ then $\psi v_2=(\varphi_1-\lambda)v_0$ for some $v_0$ and then 
        
        $$ \pmatrix{ \varphi_1  & 0\cr
                           0 & \varphi_2}\pmatrix{v_0 \cr
        -v_2}=\pmatrix{\lambda v_1+\psi v_2-\psi v_2 \cr
        -\lambda v_2}=\lambda \pmatrix{v_0 \cr
        -v_2}.$$
        Together with $(v_1,0)$ where $\varphi_1 v_1=\lambda v_1$ these span a two-dimensional eigenspace. 
     \end{prf}          
     
     The proposition shows that the divisor $D$ corresponds to the points of $S$ at which the generalized line bundle on $X$ fails to  be locally free. This data      is conserved by the Nahm flow. One may say that for a reduced curve, the singularities are part of the characteristic equation of $\phi$ and clearly conserved under the flow. For the ribbon it is the singularities of the {\it sheaf} which are conserved.   
                         
\vskip 1cm
 {Mathematical Institute,
Radcliffe Observatory Quarter,
Woodstock Road,
Oxford, OX2 6GG}
\end{document}